\documentclass[s5,11pt]{article}
\usepackage{amsmath}
\newtheorem{theorem}{Theorem}

\newtheorem{lemma}{Lemma} 

  \numberwithin{equation}{section}
\author{ Tord Sj\"odin}
\title{ A short proof of d'Alembert's Theorem. }

\begin{document}

\maketitle
\begin{abstract} A classical theorem of d'Alembert states that if a polynomial $P(x)$ with real coefficients has a non-real root $x=a+ib$, then it also has a root $x=a-ib$. We give a short and elementary inductive proof that avoids any properties of the complex conjugation operator.\end{abstract}
\paragraph{ \it AMS 2010 Subject Classsification:} Primary 12 - 01
Secondary 30 - 01
\paragraph{ \it Key words and phrases:} Polynomial, coefficient, root, complex unit, complex number \section{ Introduction.} Among the many theorems about roots of polynomials, d'Alembert's theorem on conjugate roots of polynomials with real coefficients stands out as one of the best known. The theorem and its standard proof based on the well known properties of the complex conjugation operator $z\rightarrow \overline z$ are reproduced in most elementary textbooks, see for example \cite{SRW} Ch. 4.5. A proof in terms of the isomorphism $z\rightarrow \overline z$ of the complex plane is found in \cite{BM}, Lemma, Ch.V, Section 4. We propose the following short and elementary proof using a simple inductive argument. Our notation is standard, where $\mathbf R$ denotes the real numbers, $i$ is the complex unit and $i\cdot \mathbf R$ denotes the purely imaginary numbers.
\begin{theorem} Let $P(x)=a_0+a_1x+a_2x^2+\cdots +a_nx^n$ be a polynomial with real coefficients. Assume that $P(x)$ has a non-real root $x=a+ib$, $a,b\in{\bf R}$, then $P(x)$ also has a root $x=a-ib$.
\end{theorem} 
The following lemma contains the main step in our proof of d'Alembert's theorem.
\begin{lemma} Let $P(x)=a_0+a_1x+a_2x^2+\cdots +a_nx^n$ be a polynomial with real coefficients and let $a,b$ be real numbers, then $P(a+ib)+P(a-ib)\in \mathbf {R}$ and $P(a+ib)-P(a-ib)\in i\cdot \mathbf {R}$.
\end{lemma}
Proof. Since $P(a\pm ib)=a_0+a_1(a\pm ib)+a_2(a\pm ib)^2+\cdots +a_n(a\pm ib)^n$,
 with real coefficients $a_0,a_1,a_2,\dots ,a_n$, it is sufficient to prove that
 $$ A_k=(a+ib)^k+(a-ib)^k\in  \mathbf{R}\quad  \textrm{and}\quad   B_k=(a+ib)^k-(a-ib)^k\in i\cdot \mathbf{R} \quad (\star )$$  
holds  for $0\leq k\leq n.$
 For every natural number $k$ we let $\mathcal P_k$ denote the statement that ($\star)$ is true for all real numbers $a,b$, and proceed by induction over $k$. Clearly, $\mathcal P_0$ is true since $A_0=2$ and $B_0=0$. For any natural number $k$ we get the recursion formulas
 $$A_{k+1}=a\cdot A_k+ib\cdot B_k\quad \textrm{and}\quad B_{k+1}=a\cdot B_k-ib\cdot A_k,$$
 from a straight forward calculation. Fix any natural number $k$ and assume that $\mathcal P_k$ holds. Then $A_k$ is a real number and $B_k$ is a purely imaginary number and hence $A_{k+1}$ is a real number by the recursion formula. A similar argument proves that $B_{k+1}$ is a purely imaginary number, that is $\mathcal P_{k+1}$ holds. Appealing to the Induction Axiom we conclude that $\mathcal P_k$ holds for all natural numbers $k$, which proves the lemma.\hfill $\triangle$
 \\[1em]
 {\it Proof of the theorem.} Let $P(x)$ be as in the theorem and assume that $P(a+ib)=0$, for some $a,b\in \mathbf R$. Then by the lemma, $P(a-ib)$ belongs to both $\mathbf R$ and $ i\cdot \mathbf R$. This means that $P(a-ib)=0$ and completes the proof of the theorem.\hfill $\triangle$
 \section{Historical remarks.} The french mathematician Jean d'Alembert, 1717 - 1783, was born as an illegitimate child and raised by a friend of his father. After studies in teology, law and mathematics, he became a member of the Paris Academy of Sciences in 1741 and was appointed to the French Academy in 1754, where he became its perpetual secretary in 1772. d'Alembert is perhaps best known for his work in differential equations and mechanics. He was also the first to make a serious attempt to prove the Fundamental Theorem of Algebra (FTA) in 1746, stating that every n-th degree polynomial has n complex roots. Although his proof contained gaps, it became a model for several later proofs by C. F. Gauss and others. Apart from mathematics, d'Alembert is also known as a central figure in the 18th century enlightenment and for his work with the french philosopher Denis Diderot, as editor of the scientific articles in his Encyclop\'edie. See \cite{H} and the references contained there.
   
Address: Tord Sj\"odin, Department of mathematics and mathematical statistics, University of 
Ume\aa ,  
  S-901 87 Ume\aa , Sweden;  tord.sjodin@math.umu.se

\end{document}